\documentclass[number,seceqn,dvips]{arxbj}
\usepackage{cursive}  

\aid{0} \volume{16} \issue{2} \pubyear{2010} \firstpage{343}
\lastpage{360} \doi{10.3150/09-BEJ215}

\makeatletter

\newtheorem{theo}{Theorem}[section]

\newremark{ex}[theo]{Example}

\newtheorem{lem}[theo]{Lemma}

\makeatother

\begin{document}
\begin{frontmatter}

\title{Absolute continuity for some one-dimensional processes}
\runtitle{Absolute continuity for some one-dimensional processes}

\begin{aug}
\author{\fnms{Nicolas} \snm{Fournier}\thanksref{e1}\ead[label=e1,mark,text={nicolas.fournier}]{nicolas.fournier@univ-paris12.fr}}
\and
\author{\fnms{Jacques} \snm{Printems}\thanksref{e2}\ead[label=e2,mark]{printems@univ-paris12.fr}\corref{}}
\runauthor{N. Fournier and J. Printems}
\address{Universit\'{e} Paris-Est, Laboratoire d'Analyse et de Math\'
{e}matiques Appliqu\'{e}es,
CNRS UMR 8050, Facult\'{e} des Sciences et Technologies, 61 avenue du
G\'{e}n\'{e}ral de Gaulle,
94010 Cr\'{e}teil Cedex, France.  E-mails: \printead*{e1}, \printead*{e2}}
\end{aug}

\received{\smonth{9} \syear{2008}}
\revised{\smonth{4} \syear{2009}}

%
\begin{abstract}
We introduce an elementary
method for proving the absolute continuity of
the time marginals of one-dimensional processes.
It is based on a comparison between the Fourier transform of such time
marginals with those of the one-step Euler approximation of the
underlying process.
We obtain some absolute continuity results
for stochastic differential equations with H\"{o}lder continuous coefficients.
Furthermore, we allow such coefficients to be random and to depend on the
whole path of the solution. We also show how it can be extended to
some stochastic partial differential equations and to some
L\'{e}vy-driven stochastic differential equations.
In the cases under study, the Malliavin calculus
cannot be used, because the solution in generally not Malliavin differentiable.
\end{abstract}

%
\pdfkeywords{absolute continuity, Holder coefficients, Levy processes,
random coefficients, stochastic differential equations, stochastic partial differential equations}

\end{frontmatter}

\section{Introduction}\label{intro}

In this paper, we introduce a new method for proving the absolute continuity
of the time marginals of some one-dimensional processes. The main idea
is elementary and quite rough. It is based on the explicit law of the
associated one-step Euler scheme and related to an estimate which says
that the process
and its Euler scheme remain very close to each other during one step.

As we will see, this method is quite robust and applies to many processes
for which the use of the Malliavin calculus (see Nualart \cite{n}, Malliavin
\cite{m}) is not possible because
the processes do not have Malliavin derivatives: examples of this include
SDEs with H\"{o}lder coefficients and SDEs with random coefficients.

However, we are not able, for the moment, to extend it to multidimensional
processes. The difficulty seems to be that we use some integrability
properties of some Fourier transforms which depend heavily on the
dimension.

To illustrate this method, we will consider four types of one-dimensional
processes. Let us summarize roughly the results we obtain and compare
them to
existing results.

\subsection*{Brownian SDEs with H\"{o}lder coefficients}

To introduce our method in a
simple way, we consider a process satisfying an SDE of the form
$\mathrm{d}X_t=\sigma(X_t)\,\mathrm{d}B_t + b(X_t) \,\mathrm{d}t$.
We assume that $b$ is measurable
with at most linear growth and that $\sigma$ is H\"{o}lder continuous with
exponent $\theta>1/2$. We show that $X_t$ has a density on $\{\sigma
\ne0\}$
whenever $t>0$. The proof is very short.

Such a result is probably not far from being already known. In the case
where $\sigma$ is bounded below, Aronson \cite{a}
obtains some absolute continuity results assuming
only that $\sigma$ and $b$ are measurable (together with some growth
conditions)
by analytical methods. Our result might be deduced from
\cite{a} by a localization argument, however, we did not succeed
in this direction. In any case, our proof is much simpler.

Let us observe that, to our knowledge,
all of the probabilistic papers on this topic assume at
least that $\sigma,b$ are Lipschitz continuous; see the paper by
Bouleau and Hirsch
\cite{bh} (the case where $b$ is measurable can also be treated by
using Girsanov's theorem).

Finally, let us mention that in \cite{bh}, one gets the absolute
continuity of
the law of
$X_t$ for all $t>0$ provided $\sigma(x_0)\ne0$, if $X_0=x_0$. Such a
result cannot hold in full generality for H\"{o}lder continuous coefficients:
choose $x_0>0$,
$\sigma(x)=x$ and $b(x)=-\operatorname{sign}(x) |x|^\alpha$ for some $\alpha\in(0,1)$.
Let $\tau_{\varepsilon} =\inf\{t\geq0, X_t=\varepsilon\}$ for
$\varepsilon\in\mathbb{R}_+$.
One can check, using It\^{o}'s formula,\vspace{1pt} that for $\varepsilon\in(0,x_0)$,
$\mathbb{E}[X_{t\land\tau_\varepsilon}^{1-\alpha}]
= x_0^{1-\alpha} - \mathbb{E}[\int_{0}^{t\land\tau_\varepsilon}
(\frac
{\alpha(1-\alpha)}{2}
X_s^{1-\alpha} + (1-\alpha)) \,\mathrm{d}s ] \leq x_0^{1-\alpha}-(1-\alpha)
\mathbb{E}[\tau_\varepsilon\land t]$,
whence
$\mathbb{E}[\tau_\varepsilon] \leq x_0^{1-\alpha}/(1-\alpha)$.
As a consequence, $\mathbb{E}[\tau_0] \leq x_0^{1-\alpha}/(1-\alpha
)$. But
it also holds that $X_{\tau_0+t}=0$ a.s.~for all $t\geq0$.
Thus, $\Pr[X_t=0]>0$, at
least for sufficiently large $t$.\vspace{-1pt}

\subsection*{Brownian SDEs with random coefficients depending on the
paths}

We consider here a process solving an SDE of the form
$\mathrm{d}X_t= \sigma(X_t)\kappa(t,(X_u)_{u\leq t},H_t)\,\mathrm{d}B_t+ b(t,(X_u)_{u\leq
t},H_t)\,\mathrm{d}t$ for some auxiliary adapted process $H$.
We assume some H\"{o}lder conditions on $\sigma\kappa$, some growth conditions
and that $\kappa$ is bounded below.
We prove the absolute continuity of the law of $X_t$ on the set
$\{\sigma\ne0\}$ for all $t>0$.

Observe that we do not assume that $H$ is Malliavin differentiable, which
would, of course, be needed if we wanted to use Malliavin calculus.

SDEs with random coefficients arise, for example, in finance. Indeed,
stochastic volatility models are now widely used; see, for example,
Heston \cite{h}, Fouque, Papanicolaou and Sircar
\cite{fps}.
SDEs with coefficients depending on the paths of the solutions
arise in random mechanics: if one writes an SDE satisfied by the velocity
of a particle, the coefficients will often depend on its position,
which is nothing but the integral of its velocity. One can also imagine
a particle with position $X_t$ whose diffusion and drift coefficients
depend on the distance covered by the particle
at time $t$, that is, $\sup_{[0,t]} X_s - \inf_{[0,t]} X_s$.

Here, again, the result is not far from being known: if $\sigma\kappa$
is bounded
below, one may use the result
of Gyongy \cite{g} which says that
the solution of an SDE (with random coefficients depending on the whole
paths of the solution) has the same time marginals
as the solution of an SDE with deterministic coefficients
depending only on time and position. These coefficients being measurable
and uniformly elliptic, one may then use the result of Aronson \cite{a}.
However, our method is extremely simple
and we do not have to assume that $\sigma$ is bounded below.

\subsection*{Stochastic heat equation}

We also study the heat equation
$\partial_tU = \partial_{xx}U + b(U) + \sigma(U)\dot W$ on
$\mathbb{R}_+\times[0,1]$,
with Neumann boundary
conditions, where $W$ is a space--time white noise; see Walsh \cite{w}.
We prove that $U(t,x)$ has a density on $\{\sigma\ne0\}$ for all
$t>0$ and all
$x\in[0,1]$, provided that $\sigma$ is
H\"{o}lder continuous with exponent
$\theta>1/2$ and that $b$ is measurable and has at most linear growth.

This result shows the robustness of our method: the best absolute continuity
result was due to Pardoux and Zhang \cite{pz}, who assume that $b$ and
$\sigma$
are Lipschitz continuous. Let us, however, mention that their non-degeneracy
condition is very sharp since they obtain the absolute continuity
of $U(t,x)$ for all $t>0$ and all $x\in[0,1]$, assuming only that
$\sigma(U(0,x_0))\ne0$ for some $x_0\in[0,1]$ (if $U(0,\cdot)$ is continuous).

\subsection*{L\'{e}vy-driven SDEs}

We finally consider the SDE
$\mathrm{d}X_t=\sigma(X_t)\,\mathrm{d}L_t + b(X_t)\,\mathrm{d}t$, where
$(L_t)_{t\geq0}$ is a L\'{e}vy martingale process
without Brownian part and with L\'{e}vy measure $\nu$.
Roughly, we assume that $\int_{|z|\leq\varepsilon}z^2\nu(\mathrm{d}z)\simeq
\varepsilon^{2-\lambda}$
for all $\varepsilon\in(0,1]$ and some
$\lambda\in(3/4,2)$. We obtain that the law of $X_t$ has a density on
$\{\sigma\ne0\}$ for all $t>0$, under the following assumptions:

\begin{longlist}
\item[(a)] if $\lambda\in(3/2,2)$, then $b$ is measurable and has at most
linear growth
and $\sigma$ is H\"{o}lder continuous with exponent $\theta>1/2$;

\item[(b)] if $\lambda\in[1,3/2]$, then $b$ and $\sigma$ are H\"{o}lder continuous
with exponents $\alpha>3/2-\lambda$ and $\theta>1/2$;

\item[(c)] if $\lambda\in(3/4,1)$, then $b,\sigma$ are H\"{o}lder
continuous with
exponent $\theta>3/(2\lambda)-1$.
\end{longlist}

This result appears to be the first absolute continuity result for jumping
SDEs
with non-Lipschitz coefficients.
Observe that, in some cases, we allow the drift coefficient to be only
measurable,
even when the driving L\'{e}vy process has no Brownian part. Such a result
cannot be obtained using a trick like Girsanov's theorem (because even
the law of such a L\'{e}vy process
$(L_t)_{t\in[0,1]}$ and that of $(L_t+t)_{t\in[0,1]}$
are clearly not equivalent). To our
knowledge, this gives the first absolute continuity result for L\'{e}vy-driven
SDEs with measurable drift.

Also, observe that
we allow the intensity measure of the Poissonian part to be singular:
even without Brownian part and without drift,
our result yields some absolute continuity
for L\'{e}vy-driven SDEs, even when the L\'{e}vy measure of the driving process
is completely singular. Such cases are not included in the famous works
of Bichteler and Jacod \cite{bj} or Bichteler, Gravereaux and Jacod
\cite{bgj}.
Picard \cite{p} obtained some very complete results in that direction for
SDEs with smooth coefficients. Note that Picard obtained his results
for any $\lambda\in(0,2)$: our assumption is quite heavy since
we have to restrict our study to the case where $\lambda>3/4$.

Ishikawa and Kunita \cite{ik} have obtained some
regularity results under some very simple assumptions for a different
type of jumping SDE, namely {\it canonical SDEs with jumps}; see
\cite{ik}, formula (6.1).

Let us finally mention a completely different approach developed by
Denis \cite{d}, Nourdin and Simon \cite{ns}, Bally \cite{b}, Kulik
\cite{k1,k2}
and others, where
singular L\'{e}vy measures are allowed when the drift coefficient is
sufficiently
non-constant. The case under study is truly different
since we allow the drift coefficient
to be completely degenerate.

We will frequently use the following classical lemma.

\begin{lem}\label{fourier}
For $\mu$ a non-negative finite measure on $\mathbb{R}$, we denote by
$\widehat\mu(\xi)= \int_\mathbb{R}\mathrm{e}^{\mathrm{i}\xi x}\mu(\mathrm{d}x)$ its Fourier
transform (for all $\xi\in\mathbb{R}$).
If $\int_\mathbb{R}|\widehat\mu(\xi)|^2\, \mathrm{d}\xi<\infty$,
then $\mu$ has a density with respect to the Lebesgue measure.
\end{lem}

\begin{pf}
For $n\geq1$, consider $\mu_n=\mu\star g_n$, where $g_n$ is the centered
Gaussian distribution with variance $1/n$. Then, of course,
$|\widehat\mu_n (\xi) | \leq|\widehat\mu(\xi)|$.
Furthermore, $\mu_n$ has a
density $f_n\in L^1\cap L^\infty(\mathbb{R},\mathrm{d}x)$
(for each fixed $n\geq1$), so we may
apply the Plancherel equality, which yields $\int_\mathbb
{R}f_n^2(x)\,\mathrm{d}x =
(2\curpi)^{-1}\int_\mathbb{R}|\widehat\mu_n(\xi)|^2\, \mathrm{d}\xi
\leq(2\curpi)^{-1}\int_\mathbb{R}|\widehat\mu(\xi)|^2\, \mathrm{d}\xi
=:C<\infty$. Due to the weak compactness of the balls of $L^2(\mathbb{R}
,\mathrm{d}x)$, we
may extract a subsequence $n_k$ and find a function $f\in L^2(\mathbb
{R},\mathrm{d}x)$ such
that $f_{n_k}$ goes weakly in $L^2(\mathbb{R},\mathrm{d}x)$ to $f$. But, on the
other hand,
$\mu_n(\mathrm{d}x)=f_n(x)\,\mathrm{d}x$ tends weakly (in the sense of measures)
to $\mu$.
As a consequence, $\mu$ is nothing but $f(x)\,\mathrm{d}x$.
\end{pf}

Observe here that this lemma is optimal. Indeed, the fact that
$\widehat\mu$ belongs to $L^p$ with $p>2$ does not imply that $\mu$
has a density; see counterexamples in Kahane and Salem \cite{ks}.
The following localization argument will also be of constant use.

\begin{lem}\label{loc}
For $\delta>0$, we introduce a function $f_\delta\dvtx \mathbb
{R}_+\mapsto[0,1]$,
vanishing on $[0,\delta]$, positive on $(\delta,\infty)$ and globally
Lipschitz continuous (with Lipschitz constant $1$).

Consider a probability measure $\mu$ on $\mathbb{R}$ and a function
$\sigma
\dvtx \mathbb{R}
\mapsto\mathbb{R}_+$. Assume that for each $\delta>0$, the measure
$\mu_\delta(\mathrm{d}x)=f_\delta(\sigma(x))\mu(\mathrm{d}x)$ has a density. Thus,
$\mu$ has a density on $\{x\in\mathbb{R},  \sigma(x) > 0\}$.
\end{lem}

\begin{pf}
Let $A\subset\mathbb{R}$ be a Borel set with Lebesgue measure $0$.
We have to prove that $\mu(A\cap\{\sigma>0\})=0$.
For each $\delta>0$, the measures
$\mathbf{1}_{\{\sigma(x)>\delta\}}\mu(\mathrm{d}x)$ and $\mu_\delta(\mathrm{d}x)$ are
clearly equivalent. By assumption, $\mu_\delta(A)=0$
for each $\delta>0$, whence $\mu(A\cap\{\sigma>\delta\})=0$.
Hence, $\mu(A\cap\{\sigma>0\})=\lim_{\delta\to0}
\mu(A\cap\{\sigma>\delta\})=0$.
\end{pf}

The sections of this paper are almost
independent. In Section \ref{BSDE}, we consider the case of simple Brownian
SDEs. Section \ref{BSDERC} is devoted to Brownian
SDEs with random coefficients depending on the whole path of the solution.
The stochastic heat equation is treated in Section \ref{SPDE}. Finally,
we consider some L\'{e}vy-driven SDEs in Section \ref{LSDE}.

\section{Simple Brownian SDEs}\label{BSDE}

We consider a filtered probability space $(\Omega, \mathcal{F},
(\mathcal{F}
_t)_{t\geq0},P)$
and an $(\mathcal{F}_t)_{t\geq0}$-Brownian motion $(B_t)_{t\geq0}$. For
$x\in\mathbb{R}$ and
$\sigma,b\dvtx \mathbb{R}\mapsto\mathbb{R}$,
we consider the one-dimensional SDE
%
\begin{equation}\label{sde1}
X_t = x + \displaystyle\int_0^t\sigma(X_s)\,\mathrm{d}B_s + \displaystyle\int
_0^tb(X_s)\,\mathrm{d}s.
\end{equation}

Our aim in this section is to prove the following result.

\begin{theo}\label{mt1}
Assume that $\sigma$ is H\"{o}lder continuous with exponent $\theta
\in(1/2,1]$
and that $b$ is measurable and has at most linear growth. Consider
a continuous $(\mathcal{F}_t)_{t\geq0}$-adapted solution
$(X_t)_{t\geq0}$ to (\ref{sde1}). Then, for all
$t>0$, the law of $X_t$ has a density on the set
$\{x\in\mathbb{R}, \sigma(x)\ne0\}$.
\end{theo}


Observe that the (weak or strong) existence of solutions to (\ref{sde1})
does not hold under the assumptions of Theorem \ref{mt1}.
However, at least weak existence holds if one additionally assumes
that $b$ is continuous or that $\sigma$ is bounded below; see
Karatzas and Shreve \cite{ks2}.

\begin{pf} By a scaling argument, it suffices to consider the case $t=1$.
We divide the proof into three parts.

\textit{Step} 1. For every $\varepsilon\in(0,1)$, we consider the random variable
\begin{eqnarray*}
Z_\varepsilon:= X_{1-\varepsilon}+ \int_{1-\varepsilon}^1 \sigma
(X_{1-\varepsilon})\,\mathrm{d}B_s=
X_{1-\varepsilon}+ \sigma(X_{1-\varepsilon}) (
B_{1}-B_{1-\varepsilon}).
\end{eqnarray*}
Conditioning with respect to
$\mathcal{F}_{1-\varepsilon}$, we get, for all $\xi\in\mathbb{R}$,
\begin{eqnarray*}
|\mathbb{E}[\mathrm{e}^{\mathrm{i}\xi Z_\varepsilon} \vert\mathcal
{F}_{1-\varepsilon}]| =
\bigl|\exp\bigl(\mathrm{i}\xi X_{1-\varepsilon} - \varepsilon\sigma
^2(X_{1-\varepsilon}) \xi^2/2\bigr)\bigr|
=\exp\bigl(-\varepsilon\sigma^2(X_{1-\varepsilon}) \xi^2/2\bigr).
\end{eqnarray*}

\textit{Step} 2. Using classical arguments (Doob's inequality and
Gronwall's lemma)
and the fact that $\sigma$ and $b$ have at most linear growth,
one may show that there exists a constant $C$ such that for all
$0\leq s \leq t \leq1$,
%
\begin{equation}\label{stesti}
\mathbb{E}\Bigl[\sup_{[0,1]} X_t^2\Bigr] \leq C,\qquad\mathbb
{E}
[(X_t-X_s)^2]
\leq C(t-s).
\end{equation}
Next, since $\sigma$ is H\"{o}lder continuous with index $\theta\in(1/2,1]$
and since $b$ has at most linear growth, we get,
for all $\varepsilon\in(0,1)$,
\begin{eqnarray*}
\mathbb{E}[(X_1-Z_{\varepsilon})^2] &\leq&2 \int_{1-\varepsilon}^1
\mathbb{E}
\bigl[
\bigl(\sigma(X_s)-\sigma(X_{1-\varepsilon})\bigr)^2 \bigr]\,\mathrm{d}s
+ 2 \mathbb{E}\biggl[\biggl( \int_{1-\varepsilon}^1 b(X_s)\,\mathrm{d}s
\biggr)^2\biggr] \nonumber
\\
&\leq& C \int_{1-\varepsilon}^1 \mathbb{E}
[|X_s-X_{1-\varepsilon
}|^{2\theta} ]\,\mathrm{d}s
+ 2 \varepsilon\int_{1-\varepsilon}^1 \mathbb
{E}[b^2(X_s)]\,\mathrm{d}s\nonumber
\\
&\leq& C \int_{1-\varepsilon}^1 \mathbb{E}
[|X_s-X_{1-\varepsilon}|^{2}
]^\theta \,\mathrm{d}s
+ C \varepsilon\int_{1-\varepsilon}^1 \mathbb
{E}[1+X_s^2]\,\mathrm{d}s\nonumber
\\
&\leq& C \varepsilon^{1+\theta} + C \varepsilon^2 \leq C \varepsilon
^{1+\theta},
\end{eqnarray*}
where we have used (\ref{stesti}).

\textit{Step} 3. Let $\delta>0$ be fixed,
consider the function $f_\delta$ defined in Lemma \ref{loc} and the
measure $\mu_{\delta,X_1}(\mathrm{d}x)=f_\delta(|\sigma(x)|)\mu_{X_1}(\mathrm{d}x)$, where
$\mu_{X_1}$ is the law of $X_1$. Then, for all $\xi\in\mathbb{R}$, all
$\varepsilon\in(0,1)$,
we may write
\begin{eqnarray*}
| \widehat{\mu_{\delta,X_1}} (\xi)| &=&
|\mathbb{E}[\mathrm{e}^{\mathrm{i}\xi X_1}f_\delta(|\sigma(X_1)|)]|\nonumber
\\
&\leq&|\mathbb{E}[\mathrm{e}^{\mathrm{i}\xi X_1}f_\delta(|\sigma(X_{1-\varepsilon})|)]|
+ \mathbb{E}\bigl[\bigl|f_\delta(|\sigma(X_1)|)-f_\delta(|\sigma
(X_{1-\varepsilon
})|)\bigr|\bigr]\nonumber
\\
&\leq&|\mathbb{E}[\mathrm{e}^{\mathrm{i}\xi Z_\varepsilon}f_\delta(|\sigma
(X_{1-\varepsilon})|)]|
+ |\xi| \mathbb{E}[|X_1-Z_\varepsilon|]
\\
&&{}+ \mathbb{E}\bigl[\bigl|f_\delta(|\sigma(X_1)|)-f_\delta(|\sigma
(X_{1-\varepsilon})|)\bigr|\bigr],
\end{eqnarray*}
where we used the inequality $|\mathrm{e}^{\mathrm{i}\xi x}-\mathrm{e}^{\mathrm{i}\xi z}|\leq|\xi|\cdot|x-z|$
and the fact that $f_\delta$ is bounded by $1$.
First, Step 1 implies that
\begin{eqnarray*}
|\mathbb{E}[\mathrm{e}^{\mathrm{i}\xi Z_\varepsilon} f_\delta(|\sigma
(X_{1-\varepsilon})|)]|
&\leq& \mathbb{E}[ | \mathbb{E}[\mathrm{e}^{\mathrm{i}\xi
Z_\varepsilon}f_\delta
(|\sigma(X_{1-\varepsilon})|)
\vert\mathcal{F}_{1-\varepsilon}] | ]
\\
&\leq& \mathbb{E}\bigl[ f_\delta(|\sigma(X_{1-\varepsilon})|)
\mathrm{e}^{-\varepsilon\sigma^2(X_{1-\varepsilon})\xi^2/2}\bigr]
\leq\exp(-\varepsilon\delta^2 \xi^2/2)
\end{eqnarray*}
since $f_\delta$ is bounded by $1$ and vanishes on $[0,\delta]$.
Step 2 implies that
$|\xi| \mathbb{E}[|X_1-Z_\varepsilon|] \leq C |\xi| \varepsilon
^{(1+\theta)/2}$.
Since $f_\delta$ is Lipschitz continuous and $\sigma$ is H\"{o}lder continuous
with index $\theta\in(1/2,1]$, we deduce from (\ref{stesti}) that
$\mathbb{E}[|f_\delta(|\sigma(X_1)|)-f_\delta(|\sigma
(X_{1-\varepsilon
})|)|]\leq
C\mathbb{E}[|X_1-X_{1-\varepsilon}|^\theta]\leq C \varepsilon
^{\theta/2}$.

As a conclusion, we deduce that for
all $\xi\in\mathbb{R}$ and all $\varepsilon\in(0,1)$,
\begin{eqnarray*}
| \widehat{\mu_{\delta,X_1}} (\xi)|\leq
\exp(-\varepsilon\delta^2 \xi^2/2) + C |\xi| \varepsilon
^{(1+\theta)/2}
+ C \varepsilon^{\theta/2}.
\end{eqnarray*}
For each $|\xi|\geq1$ fixed, we apply this formula with the choice
$\varepsilon:= (\log|\xi|)^2 / \xi^2 \in(0,1)$. This gives
\begin{eqnarray*}
| \widehat{\mu_{\delta,X_1}} (\xi)|\leq
\exp\bigl(- \delta^2 (\log|\xi|)^2/2\bigr) + C (\log|\xi|)^{1+\theta}/|\xi
|^{\theta}
+ C (\log|\xi|)^\theta/ |\xi|^\theta.
\end{eqnarray*}
This holding for all $|\xi|\geq1$, and $\widehat{\mu_{\delta,X_1}}$
being bounded
by $1$, we get that
$\int_\mathbb{R}| \widehat{\mu_{\delta,X_1}} (\xi)|^2 \,\mathrm{d}\xi
<\infty$
since $\theta>1/2$,
by assumption. Lemma \ref{fourier} implies that the measure $\mu
_{\delta,X_1}$
has a density for each $\delta>0$.
Lemma \ref{loc} allows us to conclude that
$\mu_{X_1}$ has a density on $\{|\sigma|>0\}$.
\end{pf}

\section{Brownian SDEs with random coefficients}\label{BSDERC}

We again start with a filtered probability space
$(\Omega, \mathcal{F}, (\mathcal{F}_t)_{t\geq0},P)$
and a $(\mathcal{F}_t)_{t\geq0}$-Brownian motion $(B_t)_{t\geq0}$.

To model the randomness of the coefficients, we consider an auxiliary
predictable process $(H_t)_{t\geq0}$, with values in some
normed space $(\mathcal{S},\| \cdot \|)$.
We then consider $\sigma\dvtx \mathbb{R}\mapsto\mathbb{R}$ and
two measurable maps $\kappa,b\dvtx \mathcal{A}\mapsto\mathbb{R}$, where
\begin{eqnarray*}
\mathcal{A}:= \{(s, (x_u)_{u\leq s}, h),   s\geq0, (x_u)_{u\geq
0}\in C(\mathbb{R}
_+,\mathbb{R}),
h \in\mathcal{S}\},
\end{eqnarray*}
and the following one-dimensional SDE:

\begin{equation}\label{sde2}
X_t=x+\int_0^t\sigma(X_s)\kappa(s,(X_u)_{u\leq s},H_s)\, \mathrm{d}B_s
+ \int_0^tb(s,(X_u)_{u\leq s},H_s) \,\mathrm{d}s.
\end{equation}

Here, again, the existence of solutions to such a
general equation
does not, of course, always hold, even under the assumptions below.
However, there are many particular cases for
which the (weak or strong) existence can be proven by
classical methods (Picard iteration, martingale problems, change of
probability, change of time, etcetera).

\begin{theo}\label{mt2}
Assume that the auxiliary process $H$ satisfies, for some $\eta>1/2$
and all $0\leq s \leq t \leq T$,
%
\begin{equation}\label{condih1}
\mathbb{E}[\|H_t\|^2]\leq C_T \quad\mbox{and} \quad
\mathbb{E}[\|H_t-H_s\|^2]\leq C_T(t-s)^\eta.
\end{equation}
Assume, also, that $\kappa\sigma$ and $b$ have at most linear growth,
that is,
for all $0\leq t \leq T$, all $(x_u)_{u\geq0}\in C(\mathbb
{R}_+,\mathbb{R})$ and all
$h\in\mathcal{S}$,
%
\begin{equation}\label{lg2}
|\sigma(x_t) \kappa(t,(x_u)_{u\leq t},h)|+
| b(t,(x_u)_{u\leq t},h)| \leq C_T\Bigl(1+ \sup_{[0,t]}|x_u|+ \|h\|\Bigr),
\end{equation}
that $\sigma$ is H\"{o}lder continuous with index $\alpha\in(1/2,1]$
and that for some $\theta_1\in(1/4,1]$, $\theta_2\in(1/2,1]$,
$\theta_3\in(1/2\eta,1]$,
all $0\leq s \leq t \leq T$, all
$(x_u)_{u\geq0}\in C(\mathbb{R}_+,\mathbb{R})$ and all $h,h'\in
\mathcal{S}$,
we have\looseness=1
\begin{eqnarray}\label{hold2}
&&|\sigma(x_t)\kappa(t,(x_u)_{u\leq t},h)
- \sigma(x_s)\kappa(s,(x_u)_{u\leq s},h')|\nonumber
\\[-8pt]\\[-8pt]
&&\quad \leq C_T \Bigl( (t-s)^{\theta_1}+\sup_{u\in[s,t]}|x_u-x_s|^{\theta
_2} +
\|h-h'\|^{\theta_3} \Bigr).\nonumber
\end{eqnarray}
Finally, assume that $\kappa$ is bounded below by some constant
$\kappa_0>0$.
Consider a continuous $(\mathcal{F}_t)_{t\geq0}$-adapted solution
$(X_t)_{t\geq0}$
to (\ref{sde2}). The law of $X_t$ then has a
density on $\{x\in\mathbb{R},   \sigma(x) \ne0\}$ whenever $t>0$.
\end{theo}

Note that (\ref{condih1}) does not imply that $H$ is a.s.~continuous:
it is
just a type of $L^2$-continuity.
Also, observe that we assume no regularity for the drift coefficient $b$.
This is not so surprising, if we consider Girsanov's theorem. However,
Girsanov's theorem might be difficult to use in such a context due
to the randomness of the coefficients (a change of probability also changes
the law of the auxiliary process).
Let us briefly illustrate (\ref{hold2}).

\begin{ex}
(a) Let $\sigma(x_s)\kappa(s,(x_u)_{u\leq s},h)=
\phi(s,x_s,\sup_{[0,s]} \varphi(x_u),h)$
with $\phi\dvtx \mathbb{R}_+\times\mathbb{R}\times\mathbb{R}\times
\mathcal{S}\mapsto\mathbb{R}$ satisfying
$|\phi(s,x,m,h)-\phi(s',x',m',h')|\leq C(|s-s'|^{\theta_1}
+|x-x'|^{\theta_2}+ |m-m'|^{\zeta} +
\|h-h'\|^{\theta_3})$ and $\varphi\dvtx \mathbb{R}\mapsto\mathbb{R}$ satisfying
$|\varphi(x)-\varphi(x')| \leq C |x-x'|^r$ with $\zeta r\geq\theta
_2$. Then,
$\sigma\kappa$ satisfies (\ref{hold2}).

(b) Let $\sigma(x_s)\kappa(s,(x_u)_{u\leq s},h)=
\phi(s,x_s,\int_0^s \varphi(x_u)\,\mathrm{d}u,h)$
with $\phi\dvtx \mathbb{R}_+ \times\mathbb{R}\times\mathbb{R}\times
\mathcal{S}\mapsto\mathbb{R}$
satisfying the
condition
$|\phi(s,x,m,h)-\phi(s',x',m',h')|\leq C(|s-s'|^{\theta
_1}+|x-x'|^{\theta_2}
+ |m-m'|^{\theta_1} + \|h-h'\|^{\theta_3})$ and with
$\varphi\dvtx \mathbb{R}\mapsto\mathbb{R}$ bounded. Then, $\sigma\kappa
$ satisfies
(\ref{hold2}).
\end{ex}

\begin{pf*}{Proof of Theorem \ref{mt2}}
The scheme of the proof is exactly the same as that of Theorem \ref{mt1}.
For the sake of simplicity, we show the result only when $t=1$.

\textit{Step} 1. For $\varepsilon\in(0,1)$, we consider the random variable
\begin{eqnarray*}
Z_\varepsilon:= X_{1-\varepsilon}+ \int_{1-\varepsilon}^1 \sigma
(X_{1-\varepsilon})\kappa\bigl(1-\varepsilon,(X_u)_{u\leq1-\varepsilon},
H_{1-\varepsilon}\bigr)\,\mathrm{d}B_s.
\end{eqnarray*}

\noindent
Conditioning with respect to
$\mathcal{F}_{1-\varepsilon}$ and using the fact that $\kappa\geq
\kappa
_0$, we get, for all
$\xi\in\mathbb{R}$,
\begin{eqnarray*}
|\mathbb{E}[\mathrm{e}^{\mathrm{i}\xi Z_\varepsilon} \vert\mathcal
{F}_{1-\varepsilon}]| &=& \bigl|
\exp\bigl(\mathrm{i}\xi X_{1-\varepsilon} - \varepsilon\sigma
^2(X_{1-\varepsilon})\kappa^2\bigl(1-\varepsilon,(X_u)_{u\leq
1-\varepsilon},
H_{1-\varepsilon}\bigr) \xi^2/2\bigr)\bigr| \nonumber
\\
&\leq& \exp\bigl(-\varepsilon\kappa_0^2\sigma^2(X_{1-\varepsilon
})\xi^2/2 \bigr).
\end{eqnarray*}

\textit{Step} 2. Using Doob's inequality, Gronwall's lemma,
(\ref{lg2}) and (\ref{condih1}), one easily shows that for all
$0\leq s \leq t \leq1$,
%
\begin{equation}\label{stesti2}
\mathbb{E}\Bigl[\sup_{[0,1]} X_t^2\Bigr] \leq C,
\qquad\mathbb{E}\Bigl[\sup_{u\in[s,t]}(X_u-X_s)^2\Bigr] \leq C(t-s).
\end{equation}
Next, using  (\ref{condih1})--(\ref{stesti2}),
we get, for all $\varepsilon\in(0,1)$,
\begin{eqnarray*}
&&\mathbb{E}[(X_1-Z_{\varepsilon})^2]
\\
&&\quad \leq 2 \int_{1-\varepsilon}^1
\mathbb{E}\bigl[
\bigl(\sigma(X_s)\kappa(s,(X_u)_{u\leq s},H_s)
-\sigma(X_{1-\varepsilon})\kappa\bigl(1-\varepsilon,(X_u)_{u\leq
1-\varepsilon},H_{1-\varepsilon}\bigr)\bigr)^2 \bigr]\,\mathrm{d}s \nonumber
\\
&&{}\qquad+ 2 \mathbb{E}\biggl[\biggl( \int_{1-\varepsilon}^1
b(s,(X_u)_{u\leq
s},H_s)\,\mathrm{d}s \biggr)^2\biggr]\nonumber
\\
&&\quad \leq C \int_{1-\varepsilon}^1 \mathbb{E}\Bigl[\bigl(s-(1-\varepsilon
)\bigr)^{2\theta_1}
+\sup_{u\in[1-\varepsilon,s]}|X_u-X_{1-\varepsilon}|^{2\theta_2} +
\|H_s-H_{1-\varepsilon}\|^{2\theta_3}\Bigr]\,\mathrm{d}s\nonumber
\\
&&{}\qquad + 2 \varepsilon\int_{1-\varepsilon}^1 \mathbb
{E}[b^2(s,(X_u)_{u\leq s},H_s)]\,\mathrm{d}s\nonumber
\\
&&\quad \leq C \varepsilon^{1+2\theta_1}
+ C \varepsilon\mathbb{E}\Bigl[\sup_{u\in[1-\varepsilon
,1]}|X_u-X_{1-\varepsilon}|^2 \Bigr]^{\theta_2}
\\
&&{}\qquad +C \varepsilon\sup_{u\in[1-\varepsilon,1]}\mathbb{E}
[\|H_u-H_{1-\varepsilon}\|^2]^{\theta_3}\nonumber
\\
&&{}\qquad + C \varepsilon\int_{1-\varepsilon}^1 \mathbb{E}\Bigl[1+\sup_{u\in
[0,s]}X_u^2+ \|H_s\|^2\Bigr] \,\mathrm{d}s \nonumber
\\
&&\quad \leq C \varepsilon^{1+2\theta_1} + C \varepsilon^{1+\theta_2}+ C
\varepsilon^{1+\eta\theta_3} + C\varepsilon^{2}
\leq C \varepsilon^{1+\theta},
\end{eqnarray*}
where $\theta:= \min(2\theta_1,\theta_2,\eta\theta_3,1) \in(1/2,1]$,
by assumption.

\textit{Step} 3. Let $\delta>0$ be fixed and consider the function
$f_\delta$
of Lemma \ref{loc} and
the measure $\mu_{\delta,X_1}(\mathrm{d}x)
=f_\delta(|\sigma(x)|)\mu_{X_1}(\mathrm{d}x)$, where
$\mu_{X_1}$ is the law of $X_1$. Then, as in the proof of Theorem \ref{mt1},
we may write, for all $\xi\in\mathbb{R}$ and all $\varepsilon\in(0,1)$,
\begin{eqnarray*}
| \widehat{\mu_{\delta,X_1}} (\xi)|
&\leq&|\mathbb{E}[\mathrm{e}^{\mathrm{i}\xi Z_\varepsilon}f_\delta(|\sigma
(X_{1-\varepsilon})|)]|
+ |\xi| \mathbb{E}[|X_1-Z_\varepsilon|]
\\
&&{}+ \mathbb{E}\bigl[ \bigl|f_\delta(|\sigma(X_1)|)-f_\delta(|\sigma
(X_{1-\varepsilon})|)\bigr|\bigr].
\end{eqnarray*}
Exactly as in the proof of Theorem \ref{mt1}, using the facts that
$\sigma$
is H\"{o}lder continuous with exponent $\alpha\in(1/2,1]$ and
that $f_\delta$ is bounded by $1$, Lipschitz continuous and vanishes
on $[0,\delta]$, we obtain from
Steps 1 and 2 that for all $\varepsilon\in(0,1)$,
\begin{eqnarray*}
| \widehat{\mu_{\delta,X_1}} (\xi)|\leq
\exp(-\varepsilon\kappa_0^2 \delta^2 \xi^2/2) + C |\xi|
\varepsilon^{(1+\theta)/2}
+ C \varepsilon^{\alpha/2}.
\end{eqnarray*}
For each $|\xi|\geq1$ fixed, we apply this formula with the choice
$\varepsilon:= (\log|\xi|)^2 / \xi^2 \in(0,1)$ and deduce, as in
the proof of
Theorem \ref{mt1}, that $\int_\mathbb{R}|
\widehat{\mu_{\delta,X_1}} (\xi)|^2 \,\mathrm{d}\xi<\infty$ because $\theta
>1/2$ and
$\alpha>1/2$.
Due to Lemma \ref{fourier}, this implies that $\mu_{\delta,X_1}$
has a density for each $\delta>0$.
Thus, $\mu_{X_1}$ has a density on $\{|\sigma|> 0\}$
thanks to Lemma \ref{loc}.
\end{pf*}

\section{Stochastic heat equation}\label{SPDE}

On a filtered probability space $(\Omega, \mathcal{F}, (\mathcal
{F}_t)_{t\geq0},P)$,
we consider an $(\mathcal{F}_t)_{t\geq0}$ space--time
white noise $W(\mathrm{d}t,\mathrm{d}x)$ on $\mathbb{R}_+\times[0,1]$,
based on $\mathrm{d}t\,\mathrm{d}x$; see Walsh \cite{w}. For two functions
$\sigma,b\dvtx \mathbb{R}\mapsto\mathbb{R}$,
we consider the stochastic heat equation
with Neumann boundary conditions
\begin{eqnarray}\label{spde}
\partial_t U(t,x)&=& \partial_{xx} U(t,x) + b(U(t,x))+ \sigma
(U(t,x))\dot W(t,x),\nonumber
\\[-8pt]\\[-8pt]
\partial_x U(t,0)&=&\partial_x U(t,1)=0,\nonumber
\end{eqnarray}
with some initial condition $U(0,x)=U_0(x)$ for some deterministic
$U_0 \in L^\infty([0,1])$.

Consider the heat kernel $G_t(x,y):=\frac{1}{\sqrt{4\curpi t}}
\sum_{n\in\mathbb{Z}}
[\mathrm{e}^{{-(y-x-2n)^2/(4t)}} + \mathrm{e}^{{-(y+x-2n)^2/(4t)}} ]$.
Following the ideas of Walsh \cite{w}, we say that a continuous
$(\mathcal{F}_t)_{t\geq0}$-adapted
process $(U(t,\break x))_{t> 0,x\in[0,1]}$ is a weak solution to
(\ref{spde}) if a.s., for all
$t>0$ and all $x\in[0,1]$,
\begin{eqnarray}\label{wspde}
U(t,x)&=&\int_0^1 G_t(x,y) U_0(y)\,\mathrm{d}y + \int_0^t\int_0^1
G_{t-s}(x,y)
b(U(s,y))\,\mathrm{d}y\,\mathrm{d}s\nonumber
\\[-8pt]\\[-8pt]
&&{}+ \int_0^t\int_0^1 G_{t-s}(x,y) \sigma(U(s,y))W(\mathrm{d}s,\mathrm{d}y).\nonumber
\end{eqnarray}

In this section, we will show the following result.

\begin{theo}\label{mt3}
Assume that $b$ is measurable and has at most linear growth, and that
$\sigma$ is
H\"{o}lder continuous with exponent $\theta\in(1/2,1]$. Consider
a continuous $(\mathcal{F}_t)_{t\geq0}$-adapted weak solution
$(U(t,x))_{t> 0,x\in[0,1]}$
to (\ref{spde}). Then, for all $x\in[0,1]$ and all $t>0$, the law of
$U(t,x)$ has
a density on $\{u \in\mathbb{R}, \sigma(u) \ne0\}$.
\end{theo}

Again, the existence of solutions is not proved under the
assumptions of Theorem \ref{mt3} alone. We mention Gatarek and Goldys
\cite{gg},
from which we obtain the weak existence of a solution
by additionally assuming that $b$ is continuous. On the other hand,
Bally, Gyongy and Pardoux \cite{bgp} have proven the existence of a solution
for a (locally) Lipschitz continuous
diffusion coefficient $\sigma$ bounded below
and a (locally) bounded
measurable drift coefficient~$b$.\looseness=1

We will use the following estimates relating to the heat kernel, which can
be found in the Appendix of Bally and Pardoux \cite{bp} and Bally,
Millet and Sanz-Sol\'{e}
\cite{bms}, Lemma~B1.
For some constants $0<c<C$,
all $\varepsilon\in(0,1)$, all $x,y \in[0,1]$ and all $0\leq\break s \leq
t \leq1$,\looseness=1
\begin{eqnarray}\label{resti1}
c \sqrt{\varepsilon} &\leq&\kappa_\varepsilon(x):=\int
_{1-\varepsilon}^1
\int_{0\lor(x-\sqrt{\varepsilon})}^{1\land(x+\sqrt{\varepsilon})}
G^2_{1-u}(x,z) \,\mathrm{d}z\,\mathrm{d}u\nonumber
\\[-8pt]\\[-8pt]
&\leq&\int_{1-\varepsilon}^1\int_0^1 G^2_{1-u}(x,z)\,
\mathrm{d}z\,\mathrm{d}u \leq C \sqrt {\varepsilon},\nonumber
\end{eqnarray}
\begin{eqnarray}\label{resti2}
&\displaystyle\int_0^t\int_0^1 \bigl(G_{t-u}(x,z)-G_{t-u}(y,z)\bigr)^2
\,\mathrm{d}z\,\mathrm{d}u \leq C |x-y|,\  &
\\\label{resti3}
&\displaystyle\int_0^s\int_0^1
\bigl(G_{t-u}(x,z)-G_{s-u}(x,z)\bigr)^2\, \mathrm{d}z\,\mathrm{d}u +
\int_s^t\int_0^1 G^2_{t-u}(x,z)\,\mathrm{d}z\,\mathrm{d}u
 \leq C |t-s|^{1/2}.\quad
\end{eqnarray}

\begin{pf*}{Proof of Theorem \ref{mt3}} We assume that $t=1$ for simplicity and we fix
$x\in[0,1]$.

\textit{Step} 1. For $\varepsilon\in(0,1)$, let
\begin{eqnarray*}
Z_\varepsilon &:=& \int_0^1 G_1(x,y) U_0(y)\,\mathrm{d}y + \int_0^{1-\varepsilon
}\int_0^1 G_{1-s}(x,y)
b(U(s,y))\,\mathrm{d}y\,\mathrm{d}s\nonumber
\\
&&{}+ \int_0^{1-\varepsilon}\int_0^1 G_{1-s}(x,y)
\sigma(U(s,y))W(\mathrm{d}s,\mathrm{d}y)
\\
&&{}+ \int_{1-\varepsilon}^1 \int_0^1 G_{1-s}(x,y) \sigma
\bigl(U(1-\varepsilon,y)\bigr)W(\mathrm{d}s,\mathrm{d}y) .
\end{eqnarray*}
As before, we observe that
\begin{eqnarray*}
|\mathbb{E}[\mathrm{e}^{\mathrm{i}\xi Z_\varepsilon}\vert\mathcal{F}_{1-\varepsilon}]|
&=& \exp\biggl(- \frac{|\xi|^2}{2}
\int_{1-\varepsilon}^1 \int_0^1 G_{1-s}^2(x,y) \sigma
^2\bigl(U(1-\varepsilon,y)\bigr)
\,\mathrm{d}y \,\mathrm{d}s \biggr)
\\
&\leq& \exp\bigl(- \kappa_\varepsilon(x) Y_\varepsilon|\xi|^2 /2\bigr),
\end{eqnarray*}
where, recalling (\ref{resti1}),
\begin{eqnarray*}
Y_\varepsilon:=\frac{1}{\kappa_\varepsilon(x)}
\int_{1-\varepsilon}^1 \int_{0\lor(x-\sqrt\varepsilon)}^{1\land
(x+\sqrt\varepsilon)} G_{1-s}^2(x,y)
\sigma^2\bigl(U(1-\varepsilon,y)\bigr)
\,\mathrm{d}y \,\mathrm{d}s.
\end{eqnarray*}

\textit{Step} 2. Using some classical computations involving
(\ref{resti1})--(\ref{resti3}), as well as the
fact that
$t,x \mapsto\int_0^1 G_t(x,y)U_0(y)\,\mathrm{d}y$ is of class $C^\infty_b$ on
$(t_0,1]\times[0,1]$ for all $t_0\in(0,1)$,
we get, for some constant $C$,
%
\begin{eqnarray}\label{stesti3.1}
&&\forall  t\in[0,1],  \forall  y\in[0,1], \qquad\mathbb{E}[U^2(t,y)]
\leq C;
 \\\label{stesti3.2}
&&\forall  s,t \in[1/2,1],  \forall  y \in[0,1], \qquad
\mathbb{E}\bigl[\bigl(U(t,y)-U(s,y)\bigr)^2\bigr]\leq C |t-s|^{1/2};
\\\label{stesti3.3}
&&\forall  t \in[1/2,1],  \forall  y,z \in[0,1], \qquad
\mathbb{E}\bigl[\bigl(U(t,y)-U(t,z)\bigr)^2\bigr]\leq C |y-z|.
\end{eqnarray}

\textit{Step} 2.1. We now prove that for all $\varepsilon\in(0,1/2)$,
\begin{eqnarray*}
\mathbb{E}\bigl[ \bigl(U(1,x) - Z_\varepsilon\bigr)^2\bigr] \leq C \varepsilon
^{(1+\theta)/2}.
\end{eqnarray*}
Since $\sigma$ is
H\"{o}lder continuous and since $b$ has at most linear growth,
using (\ref{stesti3.1}) and (\ref{stesti3.2}),
we obtain
\begin{eqnarray*}
&&\mathbb{E}\bigl[ \bigl(U(1,x) - Z_\varepsilon\bigr)^2\bigr]
\\
&&\quad \leq 2 \mathbb{E}\biggl[
\biggl(
\int_{1-\varepsilon}^1 \int_0^1 G_{1-s}(x,y)b(U(s,y))\, \mathrm{d}y\, \mathrm{d}s
\biggr)^2\biggr] \nonumber
\\
&&{}\qquad + 2 \int_{1-\varepsilon}^1 \int_0^1 G^2_{1-s}(x,y)
\mathbb{E}\bigl[\bigl(\sigma(U(s,y))-\sigma\bigl(U(1-\varepsilon
,y)\bigr)\bigr)^2\bigr]
\,\mathrm{d}y\,\mathrm{d}s\nonumber
\\
&&{}\quad \leq 2 \varepsilon\int_{1-\varepsilon}^1 \int_0^1 G^2_{1-s}(x,y)
\mathbb{E}[b^2(U(s,y))] \,\mathrm{d}y\,\mathrm{d}s \nonumber
\\
&&{}\qquad + C \int_{1-\varepsilon}^1 \int_0^1 G^2_{1-s}(x,y)
\mathbb{E}[|U(s,y)-U(1-\varepsilon,y)|^{2\theta}]\, \mathrm{d}y\,\mathrm{d}s
\nonumber
\\
&&\quad \leq C \varepsilon\int_{1-\varepsilon}^1 \int_0^1 G^2_{1-s}(x,y)
\mathbb{E}[1+U^2(s,y)]\, \mathrm{d}y\,\mathrm{d}s \nonumber
\\
&&{}\qquad + C \int_{1-\varepsilon}^1 \int_0^1 G^2_{1-s}(x,y)
\mathbb{E}[|U(s,y)-U(1-\varepsilon,y)|^{2}]^\theta\, \mathrm{d}y\,\mathrm{d}s\nonumber
\\
&&\quad \leq C \varepsilon\int_{1-\varepsilon}^1 \int_0^1 G^2_{1-s}(x,y)
\,\mathrm{d}y\,\mathrm{d}s + C \varepsilon^{\theta/2}
\int_{1-\varepsilon}^1 \int_0^1 G^2_{1-s}(x,y) \,\mathrm{d}y\,\mathrm{d}s\nonumber
\\
&&\quad \leq C \varepsilon^{3/2} + C \varepsilon^{(1+\theta)/2} \leq C
\varepsilon^{(1+\theta)/2},
\end{eqnarray*}
where, in the final inequality, we have used (\ref{resti1}).

\textit{Step} 2.2. We now check that there exists a constant
$C$ such that for all $\varepsilon\in(0,1/2)$,
\begin{eqnarray*}
A_\varepsilon:=\mathbb{E}[| \sigma^2(U(1,x)) -
Y_\varepsilon
| ]
\leq C \varepsilon^{\theta/4} .
\end{eqnarray*}
We have
\begin{eqnarray*}
A_\varepsilon & = & \frac{1}{\kappa_\varepsilon(x)} \mathbb{E}
\biggl[\biggl|
\int_{1-\varepsilon}^1 \int_{0\lor(x-\sqrt\varepsilon)}^{1\land
(x+\sqrt\varepsilon)}
G_{1-s}^2(x,y) \bigl[ \sigma^2(U(1,x))- \sigma^2\bigl(U(1-\varepsilon,y)\bigr) \bigr]\, \mathrm{d}y\,\mathrm{d}s
\biggr| \biggr]
\\
&\leq& \frac{1}{\kappa_\varepsilon(x)}
\int_{1-\varepsilon}^1 \int_{0\lor(x-\sqrt\varepsilon)}^{1\land
(x+\sqrt\varepsilon)}
G_{1-s}^2(x,y) \mathbb{E}\bigl[ \bigl|\sigma^2(U(1,x))- \sigma
^2\bigl(U(1-\varepsilon,y)\bigr)
\bigr|\bigr] \,\mathrm{d}y\,\mathrm{d}s
\\
&\leq& \sup_{y\in[x-\sqrt\varepsilon, x + \sqrt\varepsilon], }
\mathbb{E}\bigl[ \bigl|\sigma^2(U(1,x))-
\sigma^2\bigl(U(1-\varepsilon,y)\bigr) \bigr|\bigr].
\end{eqnarray*}
However, using the fact that $\sigma$ is H\"{o}lder continuous and has
at most
linear growth, using (\ref{stesti3.1})--(\ref{stesti3.3}), we deduce that for all $y\in[x-\sqrt
\varepsilon,x+\sqrt\varepsilon]$,
\begin{eqnarray*}
&&\mathbb{E}\bigl[\bigl|\sigma^2(U(1,x))-\sigma^2\bigl(U(1-\varepsilon,y)\bigr)\bigr|
\bigr]
\\
&&\quad \leq
\mathbb{E}\bigl[\bigl|\sigma(U(1,x))-\sigma\bigl(U(1-\varepsilon,y)\bigr)\bigr|^2
\bigr]^{1/2}
\mathbb{E}\bigl[\bigl|\sigma(U(1,x))+\sigma
\bigl(U(1-\varepsilon
,y)\bigr)\bigr|^2 \bigr]^{1/2}
\nonumber
\\
&&\quad \leq C \mathbb{E}[|U(1,x)-U(1-\varepsilon,y)|^{2\theta}
]^{1/2}\nonumber
\\
&&\quad \leq C \mathbb{E}[|U(1,x)-U(1-\varepsilon,y)|^{2}
]^{\theta
/2}\leq
C(\varepsilon^{1/2}+ |x-y|)^{\theta/2} \leq C \varepsilon^{\theta/4},
\end{eqnarray*}
which concludes the step.

\textit{Step} 3. Denote by $\mu_{U(1,x)}$ the
law of $U(1,x)$. For $\delta>0$, consider $f_\delta$ as in Lemma \ref
{loc} and
set $\mu_{\delta,U(1,x)}(\mathrm{d}u)=f_\delta(\sigma^2(u))\mu_{U(1,x)}(\mathrm{d}u)$.
For all $\xi\in\mathbb{R}$ and all $\varepsilon\in(0,1/2)$, we may write,
as in the proof
of Theorem \ref{mt1},
\begin{eqnarray*}
|\widehat{\mu_{\delta,U(1,x)}} (\xi)| &=& |\mathbb{E}[\mathrm{e}^{\mathrm{i}\xi U(1,x)}
f_\delta(\sigma^2((U(1,x)))]| \nonumber
\\
&\leq& |\mathbb{E}[\mathrm{e}^{\mathrm{i}\xi Z_\varepsilon}
f_\delta( Y_\varepsilon) ]|
+ |\xi| \mathbb{E}[|U(1,x) - Z_\varepsilon|] +
\mathbb{E}[ |f_\delta(\sigma^2(U(1,x)))
- f_\delta(Y_\varepsilon)|].
\end{eqnarray*}
Using Steps 1, 2.1 and 2.2, observing that $Y_\varepsilon$ is
$\mathcal{F}
_{1-\varepsilon}$-measurable
and recalling that $f_\delta$
is bounded by $1$ and vanishes on $[0,\delta]$, we get
\begin{eqnarray*}
|\widehat{\mu_{\delta,U(1,x)}} (\xi)|
\leq \mathrm{e}^{- \kappa_\varepsilon(x) \delta\xi^2/2} + C|\xi|
\varepsilon^{(1+\theta)/4}
+C \varepsilon^{\theta/4} \leq \mathrm{e}^{- c\delta\sqrt\varepsilon\xi
^2/2} + C|\xi| \varepsilon^{(1+\theta)/4}
+C \varepsilon^{\theta/4},
\end{eqnarray*}
using (\ref{resti1}) for the last inequality.
For each $|\xi|\geq1$, we choose $\varepsilon:=(\log|\xi|)^4/\xi
^4 \in(0,1/2)$
and get
\begin{eqnarray*}
|\widehat{\mu_{\delta, U(1,x)}} (\xi)|
 \leq\exp\bigl(- c \delta(\log|\xi|)^2/2 \bigr)
+ C (\log|\xi|)^{1+\theta}/|\xi|^\theta+ C (\log|\xi|)^\theta
/|\xi|^\theta.
\end{eqnarray*}
This holding for all $|\xi|\geq1$ and $|\widehat{\mu_{\delta
,U(1,x)}} (\xi)|$
being bounded by $1$, we conclude, since $\theta>1/2$, that
$\int_\mathbb{R}|\widehat{\mu_{\delta,U(1,x)}} (\xi)|^2\,\mathrm{d}\xi
<\infty$.
Lemma \ref{fourier}
ensures that the law of $\mu_{\delta,U(1,x)}$ has a density
for each $\delta>0$. We conclude,
using Lemma \ref{loc}, that $\mu_{U(1,x)}$ has a density on $\{\sigma
^2 > 0\}$.
\end{pf*}

\section{L\'{e}vy-driven SDEs}\label{LSDE}

We conclude this paper by considering L\'{e}vy-driven SDEs.
For simplicity, we restrict our study to the case of deterministic coefficients
depending only on the position of the process. The result below extends
without
difficulty, as in the Brownian case, to SDEs with random coefficients
depending on the whole paths, under some adequate conditions.

We thus
consider a filtered probability space
$(\Omega, \mathcal{F}, (\mathcal{F}_t)_{t\geq0},P)$ and
a square-integrable compensated
$(\mathcal{F}_t)_{t\geq0}$-L\'{e}vy process $(L_t)_{t\geq0}$ without
drift, without
Brownian part and with L\'{e}vy measure $\nu$. Such a process is entirely
characterized by its Fourier transform:
\begin{eqnarray*}
\mathbb{E}[\exp(\mathrm{i} \xi L_t) ]=
\exp\biggl(-t \int_{\mathbb{R}_*} (1-\mathrm{e}^{\mathrm{i}\xi z}+\mathrm{i}\xi z) \nu
(\mathrm{d}z)\biggr).
\end{eqnarray*}
For $\sigma,b\dvtx\mathbb{R}\mapsto\mathbb{R}$, we consider the
one-dimensional SDE
%
\begin{equation}\label{sde4}
X_t = x + \int_0^t\sigma(X_{s-})\,\mathrm{d}L_s +
\int_0^tb(X_s)\,\mathrm{d}s.
\end{equation}

Our aim in this section is to prove the following result.

\begin{theo}\label{mt4}
Assume that $\int_{\mathbb{R}_*} z^2\nu(\mathrm{d}z)<\infty$
and that for some $\lambda\in(3/4,2)$, $c>0$, $\xi_0\geq0$,
%
\begin{equation}\label{tasoeur1}
\forall  |\xi|\geq\xi_0\qquad   \int_{\mathbb{R}_*} \bigl(1-\cos(\xi
z)\bigr)\nu(\mathrm{d}z)
\geq c |\xi|^\lambda
\end{equation}
and for some $\gamma\in[1,2]$ (with, necessarily, $\gamma\geq
\lambda$),
%
\begin{equation}\label{tasoeur2}
\int_{\mathbb{R}_*} |z|^\gamma\nu(\mathrm{d}z) <\infty.
\end{equation}
We also assume that $b$ is measurable with at most linear growth and
that $\sigma$ is H\"{o}lder continuous with exponent
$\theta\in(3\gamma/(2\lambda) -1,1]$. If $\lambda\in(3/4,3/2)$, we
additionally suppose that $b$ is H\"{o}lder continuous with index
$\alpha\in(3\gamma/(2\lambda)-\gamma,1]$.

Let $(X_t)_{t\geq0}$ be a cadlag $(\mathcal{F}_t)_{t\geq0}$-adapted solution
to (\ref{sde4}). Then, for all $t>0$, the law of $X_t$ has a density
on the set $\{x\in\mathbb{R},  \sigma(x)\ne0\}$.
\end{theo}

Here, again, the (weak or strong) existence of solutions to (\ref{sde4})
probably does not hold under the assumptions of Theorem \ref{mt4} alone.
See Jacod \cite{j} for many existence results.

Let us comment on this result.

(a) Observe that (\ref{tasoeur2}) implies
$\int_{\mathbb{R}_*} (1-\cos(\xi z))\nu(\mathrm{d}z) \leq C |\xi|^\gamma$
so that under
(\ref{tasoeur1}),
(\ref{tasoeur2}) can hold only for some $\gamma\geq\lambda$.

Indeed, since $0 \leq1-\cos x \leq2(x^2\land1)$, we may write
$\int_{\mathbb{R}_*} (1-\cos(\xi z))\nu(\mathrm{d}z)
\leq2\int_{|z|\leq1/|\xi|} \xi^2\times z^2 \nu(\mathrm{d}z) +
2 \int_{|z|\geq1/|\xi|} \nu(\mathrm{d}z) \leq2 \xi^2 \int_{|z|\leq1/|\xi
|} |z|^{\gamma}
|\xi|^{\gamma-2} \nu(\mathrm{d}z) +\break 2 \int_{|z|\geq1/|\xi|} |z|^\gamma
|\xi|^\gamma
\nu(\mathrm{d}z) \leq2 |\xi|^\gamma\times \int_{\mathbb{R}_*} |z|^\gamma\nu(\mathrm{d}z)$.

(b) Using a standard localization procedure, one may easily eliminate
large jumps, that is, replace the assumptions
$\int_{\mathbb{R}_*}(|z|^2+|z|^\gamma)\nu(\mathrm{d}z)<\infty$ by
$\int_{\mathbb{R}_*} \min(1,|z|^\gamma)\nu(\mathrm{d}z)<\infty$.

(c) If (\ref{tasoeur1}) holds with $\lambda>3/2$, we assume
no regularity on the drift coefficient $b$.
Observe, here, that no trick using Girsanov's theorem may allow us to
remove the
drift: there is a clear difference in nature between the paths
of a L\'{e}vy process without Brownian part with and without drift.

(d) Assume that $\nu$ satisfies $\int_{\mathbb{R}_*} z^2 \nu
(\mathrm{d}z)<\infty$
and that the following
property holds for some $\lambda\in(3/4,2)$:
there exist $0<c_0<c_1$ such that for all $\varepsilon\in(0,1]$,
%
\begin{equation}\label{dix}
c_0 \varepsilon^{2-\lambda} \leq\int_{|z|\leq\varepsilon} z^2 \nu
(\mathrm{d}z) \leq c_1 \varepsilon^{2-\lambda}.
\end{equation}
Then, (\ref{tasoeur1}) holds and (\ref{tasoeur2})
holds with any $\gamma\in(\lambda,2]$. Indeed, since $1-\cos x \geq x^2/2$
for $x\in[0,1]$, we get, for $|\xi|> 1$, $\int_{\mathbb{R}_*}
(1-\cos(\xi
z)) \nu(\mathrm{d}z)
\geq(\xi^2/4) \int_{|z|\leq1/|\xi|} z^2 \nu(\mathrm{d}z) \geq c_0 |\xi
|^\lambda/4$,
whence (\ref{tasoeur1}). Next, let $\gamma\in(\lambda,2)$ be fixed.
To show that (\ref{tasoeur2}) holds, it clearly suffices to prove that
$\int_{|z|<1} |z|^\gamma\nu(\mathrm{d}z) <\infty$. Let us, for example, show that
$\int_0^1 z^\gamma\nu(\mathrm{d}z) <\infty$. Using an integration by parts,
one easily gets $\int_0^1 z^\gamma\nu(\mathrm{d}z) = \int_0^1 z^{\gamma
-2}z^2\nu(\mathrm{d}z)
= \int_0^1 (2-\gamma)z^{\gamma-3} [\int_0^z y^2\nu(\mathrm{d}y)]\, \mathrm{d}z
\leq(2-\gamma) c_1 \int_0^1 z^{\gamma-3}z^{2-\lambda} \,\mathrm{d}z <\infty$ since
$\gamma-\lambda>0$.

Thus, our result holds in the
following situations:

\begin{itemize}[$\bullet$]
\item $\lambda>3/2$, $\sigma$ is H\"{o}lder continuous with exponent
$\theta>1/2$;

\item $\lambda\in[1,3/2]$, $\sigma$ is H\"{o}lder continuous
with index
$\theta>1/2$, $b$ is H\"{o}lder continuous with exponent $\alpha
>3/2-\lambda$;

\item $\lambda\in(3/4,1]$, $\sigma$ and $b$ are H\"{o}lder
continuous with
exponent $\theta>3/(2\lambda)-1$.
\end{itemize}

(e) For example,
$\nu(\mathrm{d}z)=z^{-1-\lambda}\mathbf{1}_{[0,1]}(z)\,\mathrm{d}z$ satisfies (\ref{dix}),
as well as $\nu(\mathrm{d}z)=\sum_{n\geq1} n^{\lambda-1} \delta_{1/n}$, or, more
generally, $\nu(\mathrm{d}z)=\sum_{n\geq1} n^{\lambda\alpha-1} \delta
_{n^{-\alpha}}$ with
$\alpha>0$.

(f) Our assumption that $\lambda>3/4$ might seem strange.
However, our method does not seem to work for smaller values of
$\lambda$, even if $\sigma,b$ are Lipschitz continuous.

As noted by the anonymous referee,
however, it is possible to obtain some results for
$\lambda\in(1/2,3/4]$ if there is no drift part ($b\equiv0$).

\begin{pf*}{Proof of Theorem \ref{mt4}} By scaling, it suffices to consider the case $t=1$.
We will often write the L\'{e}vy process as
\begin{eqnarray*}
L_t=\int_0^t \int_{\mathbb{R}_{*}} z
\tilde{N}(\mathrm{d}s,\mathrm{d}z),
\end{eqnarray*}
where $\tilde{N}(\mathrm{d}s,\mathrm{d}z)$ is a compensated
Poisson measure on $\mathbb{R}_+ \times\mathbb{R}_*$ with intensity
measure $\mathrm{d}s\nu(\mathrm{d}z)$.
Thus, (\ref{sde4}) can be rewritten as
%
\begin{equation}\label{levypoisson}
X_t= x + \int_0^t \int_{\mathbb{R}_*} \sigma(X_{s-}) z
\tilde{N}(\mathrm{d}s,\mathrm{d}z) +
\int_0^tb(X_s)\,\mathrm{d}s.
\end{equation}

\textit{Step} 1. For $\varepsilon\in(0,1)$, we consider the random variable
\begin{eqnarray*}
Z_\varepsilon:= X_{1-\varepsilon}+ \int_{1-\varepsilon}^1 \sigma
(X_{1-\varepsilon})\,\mathrm{d}L_s + \int_{1-\varepsilon}^1 b(X_{1-\varepsilon})
\,\mathrm{d}s.
\end{eqnarray*}
For $\delta>0$, consider the function $f_\delta$ of Lemma
\ref{loc}. Recall that $f_\delta$ is bounded and vanishes on
$[0,\delta]$.
Conditioning with respect to
$\mathcal{F}_{1-\varepsilon}$ and using (\ref{tasoeur1}), we get,
for all $|\xi|\geq\xi_0/\delta$,
\begin{eqnarray*}
&&|\mathbb{E}[\mathrm{e}^{\mathrm{i}\xi Z_\varepsilon} f_\delta(|\sigma
(X_{1-\varepsilon})|)
\vert\mathcal{F}_{1-\varepsilon}]|\nonumber
\\
&&\quad = f_\delta(|\sigma(X_{1-\varepsilon})|)
\biggl|\exp\biggl(\mathrm{i}\xi X_{1-\varepsilon}+\mathrm{i}\xi\varepsilon
b(X_{1-\varepsilon})
\\
&&{}\qquad\hspace{80pt} - \varepsilon\int_{\mathbb{R}_*}
\bigl(1-\mathrm{e}^{\mathrm{i}\xi\sigma(X_{1-\varepsilon})z} + \mathrm{i} \xi\sigma
(X_{1-\varepsilon})z \bigr)\nu(\mathrm{d}z)\biggr)\biggr|\nonumber
\\
&&\quad = f_\delta(|\sigma(X_{1-\varepsilon})|)
\exp\biggl(- \varepsilon\int_{\mathbb{R}_*} \bigl(1-\cos(\xi\sigma
(X_{1-\varepsilon})z)\bigr)\nu(\mathrm{d}z)\biggr)\nonumber
\\
&&\quad \leq f_\delta(|\sigma(X_{1-\varepsilon})|)\exp(-c \varepsilon
\delta^{\lambda}
|\xi|^\lambda) \leq\exp(-c \varepsilon\delta^{\lambda}
|\xi|^\lambda).
\end{eqnarray*}
We have used the fact that $f_\delta$ is bounded by $1$ and vanishes
on $[0,\delta]$
to obtain the two last inequalities.

\textit{Step} 2. Recall that $\sigma$ and $b$ are H\"{o}lder continuous with
exponent
$\theta\in(0,1]$ and $\alpha\in[0,1]$ (when there is no regularity
assumption on $b$, we say that
it is H\"{o}lder with exponent $0$). The goal of this step is to show
that for all $\varepsilon\in(0,1)$,
%
\begin{equation}\label{goal}
I_\varepsilon:=\mathbb{E}[|X_1 -Z_\varepsilon|^\gamma] \leq C
\varepsilon
^{1+\theta} + C \varepsilon^{\gamma+ \alpha}
\leq C \varepsilon^{1+ \zeta},
\end{equation}
where $\zeta:=\min(\theta, \gamma+ \alpha-1) \in(3\gamma
/2\lambda-1,1]$,
by assumption.
We first show that for all $0\leq s \leq t \leq1$,
%
\begin{equation}\label{stesti4}
\mathbb{E}\Bigl[\sup_{[0,1]} |X_s|^\gamma\Bigr] \leq C, \qquad
\mathbb{E}
[|X_t-X_s|^\gamma]
\leq C |t-s|.
\end{equation}
First, using (\ref{levypoisson}),
the Burkholder--Davies--Gundy inequality (see Dellacherie and Meyer
\cite{dm}), the subadditivity of $x\mapsto x^{\gamma/2}$, the H\"{o}lder
inequality, (\ref{tasoeur2}) and the fact that $b$,
$\sigma$ have at most linear growth, we obtain, for all $t\in[0,1]$,
\begin{eqnarray*}
&&\mathbb{E}\Bigl[\sup_{u\in[0,t]} |X_u|^\gamma\Bigr]
\\[3pt]
&&\quad \leq
C|x|^\gamma+
C \mathbb{E}\biggl[\sup_{u\in[0,t]}\biggl| \int_0^u \int_{\mathbb
{R}_*} \sigma
(X_{s-})z
\tilde{N}(\mathrm{d}s,\mathrm{d}z)
\biggr|^{\gamma}\biggr]
+ C \mathbb{E}\biggl[ \biggl( \int_0^t|b(X_s)|\,\mathrm{d}s
\biggr)^\gamma\biggr] \nonumber
\\[3pt]
&&\quad \leq C|x|^\gamma+
C \mathbb{E}\biggl[\biggl(\int_0^t\int_{\mathbb{R}_*}
|\sigma(X_{s-})z|^2 N(\mathrm{d}s,\mathrm{d}z)
\biggr)^{\gamma/2}\biggr]
+ C \mathbb{E}\biggl[ \biggl( \int_0^t|b(X_s)|\,\mathrm{d}s
\biggr)^\gamma\biggr] \nonumber
\\[3pt]
&&\quad \leq C|x|^\gamma+
C \mathbb{E}\biggl[\int_0^t\int_{\mathbb{R}_*}
|\sigma(X_{s-})z|^\gamma N(\mathrm{d}s,\mathrm{d}z)\biggr]
+ C t^{\gamma-1} \mathbb{E}\biggl[\int
_0^t|b(X_s)|^\gamma \,\mathrm{d}s \biggr] \nonumber
\\[3pt]
&&\quad \leq C|x|^\gamma+
C \int_0^t\int_{\mathbb{R}_*} \mathbb{E}[|\sigma
(X_{s-})|^\gamma] |z|^\gamma\nu(\mathrm{d}z)\,\mathrm{d}s
+ C t^{\gamma-1} \int_0^t\mathbb{E}[|b(X_s)|^\gamma]
\,\mathrm{d}s\nonumber
\\[3pt]
&&\quad \leq C|x|^\gamma+C \int_0^t\mathbb{E}[1+
|X_s|^\gamma] \,\mathrm{d}s
\end{eqnarray*}
and Gronwall's lemma allows us to conclude that
$\mathbb{E}[\sup_{[0,1]} |X_s|^\gamma] \leq C$. The same arguments
ensure that
for $0\leq s \leq t \leq1$,
$\mathbb{E}[|X_t-X_s|^\gamma]\leq C \int_s^t \mathbb{E}[1+
|X_u|^\gamma] \,\mathrm{d}u$,
whence the second inequality of (\ref{stesti4}). We may now check
(\ref{goal}).
Using similar arguments and the H\"{o}lder continuity assumptions,
we obtain
\begin{eqnarray*}
I_\varepsilon
&\leq& C \mathbb{E}\biggl[\biggl(\int_{1-\varepsilon}^1 \int_{\mathbb{R}_*}
\bigl|\bigl(\sigma(X_{s-})-\sigma(X_{1-\varepsilon})\bigr)z\bigr|^2 N(\mathrm{d}s,\,\mathrm{d}z)
\biggr)^{\gamma/2}\biggr]
\\
&&{}+ C \mathbb{E}\biggl[ \biggl( \int_{1-\varepsilon}^1
|b(X_s)-b(X_{1-\varepsilon})|\,\mathrm{d}s \biggr)^\gamma\biggr]\nonumber
\\
&\leq& C \int_{1-\varepsilon}^1 \mathbb{E}[|\sigma
(X_{s-})-\sigma
(X_{1-\varepsilon})|^\gamma]\,\mathrm{d}s
+ C \varepsilon^{\gamma-1} \int_{1-\varepsilon}^1
\mathbb{E}[|b(X_{s-})-b(X_{1-\varepsilon})|^\gamma
]\,\mathrm{d}s\nonumber
\\
&\leq& C \int_{1-\varepsilon}^1 \mathbb{E}
[|X_s-X_{1-\varepsilon
}|^{\gamma\theta}]\,\mathrm{d}s
+ C \varepsilon^{\gamma-1} \int_{1-\varepsilon}^1 \mathbb{E}
[|X_s-X_{1-\varepsilon}|^{\alpha\gamma}]
\,\mathrm{d}s\nonumber
\\
&\leq& C \int_{1-\varepsilon}^1 \mathbb{E}
[|X_s-X_{1-\varepsilon
}|^\gamma]^\theta \,\mathrm{d}s
+ C \varepsilon^{\gamma-1} \int_{1-\varepsilon}^1 \mathbb{E}
[|X_s-X_{1-\varepsilon}|^\gamma]^\alpha\, \mathrm{d}s
\\
&\leq& C \varepsilon^{1+\theta}+ C \varepsilon^{\gamma+\alpha},
\end{eqnarray*}
where, in the final inequality, we have used (\ref{stesti4}).

\textit{Step} 3. Let $\delta>0$ be fixed and
consider the measure $\mu_{\delta,X_1}(\mathrm{d}x)
=f_\delta(|\sigma(x)|) \mu_{X_1}(\mathrm{d}x)$, where
$\mu_{X_1}$ is the law of $X_1$. Then, as before,
for all $\xi\in\mathbb{R}$ and all $\varepsilon\in(0,1)$, we may write
\begin{eqnarray*}
| \widehat{\mu_{\delta,X_1}} (\xi)|
&\leq& |\mathbb{E}[\mathrm{e}^{\mathrm{i}\xi Z_\varepsilon}f_\delta(|\sigma
(X_{1-\varepsilon})|)]|
+ |\xi| \mathbb{E}[|X_1-Z_\varepsilon|]
+ \mathbb{E}\bigl[ \bigl|f_\delta(|\sigma(X_1)|)-f_\delta(|\sigma
(X_{1-\varepsilon})|)\bigr|\bigr].
\end{eqnarray*}
Using the H\"{o}lder continuity of $\sigma$ and (\ref{stesti4}), one easily
gets (recall that $0< \theta\leq1 \leq\gamma$, by assumption)
$\mathbb{E}[ |f_\delta(|\sigma(X_1)|)-f_\delta(|\sigma
(X_{1-\varepsilon})|)|]
\leq C \mathbb{E}[|X_1-X_{1-\varepsilon}|^\theta] \leq C \varepsilon
^{\theta/\gamma}$.
Next, using Steps~1 and~2, we obtain,
for all $\varepsilon\in(0,1)$ and all $|\xi| \geq\xi_0/\delta$,
\begin{eqnarray*}
| \widehat{\mu_{\delta,X_1}} (\xi)|\leq
\exp(- c \delta^{\lambda} \varepsilon|\xi|^\lambda)
+ C |\xi| \varepsilon^{(1+\zeta)/\gamma} + C \varepsilon^{\theta
/\gamma}.
\end{eqnarray*}
For each $|\xi|\geq\xi_1\lor(\xi_0/\delta)$,
we choose $\varepsilon:= (\log|\xi|)^2 / |\xi|^\lambda\in(0,1)$
(this holds if $\xi_1$ is large enough). This gives
\begin{eqnarray*}
| \widehat{\mu_{\delta,X_1}} (\xi)|&\leq&
\exp(- c \delta^{\lambda} (\log|\xi|)^2 )
+ C (\log|\xi|)^{2(1+\zeta)/\gamma}/|\xi|^{\lambda(1+\zeta
)/\gamma-1}
\\
&&{}+ C (\log|\xi|)^{2\theta/\gamma} / |\xi|^{\lambda\theta/\gamma}.
\end{eqnarray*}
This holding for all $|\xi|\geq\xi_1\lor(\xi_0/\delta)$ and
$\widehat{\mu_{\delta,X_1}}$ being bounded by $1$, we get that\break
$\int_\mathbb{R}| \widehat{\mu_{\delta,X_1}} (\xi)|^2\, \mathrm{d}\xi
<\infty$.
Indeed, $\lambda(1+\zeta)/\gamma-1 >1/2$ (because
$\zeta>3\gamma/2\lambda-1$) and
$\lambda\theta/\gamma> 1/2$ (because $\theta>3\gamma/2\lambda-1$
and $\lambda\leq\gamma$).
Lemma \ref{fourier} implies that the measure $\mu_{\delta,X_1}$
has a density (for $\delta>0$ fixed) and we conclude using Lemma
\ref{loc} that
$\mu_{X_1}$ has a density on $\{|\sigma|> 0\}$.
\end{pf*}


\printhistory

\end{document}